\documentclass{amsart}
\begin{document}

\vfuzz2pt 
\hfuzz2pt 
\newtheorem{thm}{Theorem}[section]
\newtheorem{corollary}[thm]{Corollary}
\newtheorem{lemma}[thm]{Lemma}
\newtheorem{proposition}[thm]{Proposition}
\newtheorem{defn}[thm]{Definition}
\newtheorem{remark}[thm]{Remark}
\newtheorem{example}[thm]{Example}
\newtheorem{fact}[thm]{Fact}
\
\newcommand{\norm}[1]{\left\Vert#1\right\Vert}
\newcommand{\abs}[1]{\left\vert#1\right\vert}
\newcommand{\set}[1]{\left\{#1\right\}}
\newcommand{\Real}{\mathbb R}
\newcommand{\eps}{\varepsilon}
\newcommand{\To}{\longrightarrow}
\newcommand{\BX}{\mathbf{B}(X)}
\newcommand{\A}{\mathcal{A}}
\newcommand{\onabla}{\overline{\nabla}}
\newcommand{\hnabla}{\hat{\nabla}}


\def\proof{\medskip Proof.\ }
\font\lasek=lasy10 \chardef\kwadrat="32 
\def\kwadracik{{\lasek\kwadrat}}
\def\koniec{\hfill\lower 2pt\hbox{\kwadracik}\medskip}

\newcommand*{\C}{\mathbf{C}}
\newcommand*{\R}{\mathbf{R}}
\newcommand*{\Z}{\mathbf {Z}}

\def\sb{f:M\longrightarrow \C ^n}
\def\det{\hbox{\rm det}\, }
\def\detc{\hbox{\rm det }_{\C}}
\def\i{\hbox{\rm i}}
\def\tr{\hbox{\rm tr}\, }
\def\rk{\hbox{\rm rk}\,}
\def\vol{\hbox{\rm vol}\,}
\def\Im {\hbox{\rm Im}\, }
\def\Re{\hbox{\rm Re}\, }
\def\interior{\hbox{\rm int}\, }
\def\e{\hbox{\rm e}}
\def\pu{\partial _u}
\def\pv{\partial _v}
\def\pui{\partial _{u_i}}
\def\puj{\partial _{u_j}}
\def\puk{\partial {u_k}}
\def\div{\hbox{\rm div}\,}
\def\Ric{\hbox{\rm Ric}\,}
\def\r#1{(\ref{#1})}
\def\ker{\hbox{\rm ker}\,}
\def\im{\hbox{\rm im}\, }
\def\I{\hbox{\rm I}\,}
\def\id{\hbox{\rm id}\,}
\def\exp{\hbox{{\rm exp}^{\tilde\nabla}}\.}
\def\cka{{\mathcal C}^{k,a}}
\def\ckplusja{{\mathcal C}^{k+1,a}}
\def\cja{{\mathcal C}^{1,a}}
\def\cda{{\mathcal C}^{2,a}}
\def\cta{{\mathcal C}^{3,a}}
\def\c0a{{\mathcal C}^{0,a}}
\def\f0{{\mathcal F}^{0}}
\def\fnj{{\mathcal F}^{n-1}}
\def\fn{{\mathcal F}^{n}}
\def\fnd{{\mathcal F}^{n-2}}
\def\Hn{{\mathcal H}^n}
\def\Hnj{{\mathcal H}^{n-1}}
\def\emb{\mathcal C^{\infty}_{emb}(M,N)}
\def\M{\mathcal M}
\def\Ef{\mathcal E _f}
\def\Eg{\mathcal E _g}
\def\Nf{\mathcal N _f}
\def\Ng{\mathcal N _g}
\def\Tf{\mathcal T _f}
\def\Tg{\mathcal T _g}
\def\diff{{\mathcal Diff}^{\infty}(M)}
\def\embM{\mathcal C^{\infty}_{emb}(M,M)}
\def\U1f{{\mathcal U}^1 _f}
\def\Uf{{\mathcal U} _f}
\def\Ug{{\mathcal U} _g}
\def\[f]{{\mathcal U}^1 _{[f]}}
\def\hnu{\hat\nu}
\def\gnu{\nu_g}

\def\onabla{\overline\nabla}
\def\hnabla{\hat\nabla}
\def\oR{\overline R}
\def\oRic{\overline {\Ric}}
\def\hR{\hat R}

\title{Completness in affine and statistical geometry}
\author{Barbara Opozda}



\address{Faculty of Mathematics and Computer Science UJ,
ul. \L ojasiewicza  6, 30-348 Cracow, Poland}

\email{barbara.opozda@im.uj.edu.pl} \maketitle

\begin{abstract}
We begin the study of completeness of affine connections, especially
those on statistical  manifolds or on affine hypersurfaces. We
collect basic facts, prove new theorems and provide examples with
remarkable properties.

\end{abstract}
\section{Introduction}

 Although affine  hypersurfaces have been
studied for more than 100 years, the completeness of  affine
connections naturally appearing  on such hypersurfaces has not been
considered  at all. Up till now the affine completness of an affine
hypersurface has always been meant as the completeness of the affine
metric. There are a few very famous theorems in this respect in the
case of Blaschke hypersurfaces, like theorems of Blaschke, Calabi,
Cheng-Yau, see e.g. \cite{C}, \cite{CY}, \cite{LSZ}. As concerns the
completeness of the induced affine connection (which is usually
non-metrizable) it was only noticed by Nomizu, see e.g. \cite{NS},
that the two notions of completeness, that is, the completeness of
the  affine metric and the one of the induced connection,  are
independent in general.  In this paper we prove some positive and
negative results on completeness of the induced connections on
hypersurfaces. In particular, we prove that on a centroaffine
ovaloid the induced connection is complete. Of course,  connections
on compact manifolds do not have to be complete. We also observe
that on an affine hypersurface with parallel cubic form the induced
connections is not complete unless the induced structure is trivial.
In fact, we prove this result in a more general setting, namely for
statistical structures on abstract manifolds and with the assumption
weaker than that about parallel cubic form.  The class of
hypersurfaces with parallel cubic form is rich of examples and
important in affine differential geometry. Blaschke hyperurfaces
satisfying this condition were classified in \cite{HLV}.

Very little is also known  about complete statistical connections.
The study of complete statistical connections on abstract manifolds
was initiated by Noguchi in \cite{No}. He gave a procedure of
producing complete statistical connections on complete Riemannian
manifolds, see Theorem \ref{Noguchi}. In this paper we prove new
results on completeness of statistical connections. For instance, we
prove that on a topological sphere a Ricci-non-degenerate
projectively flat statistical connection  is complete.

We also provide some facts dealing with completeness of affine
connections in general. It turns out that using a cubic form defined
by  a given affine connection and some additionally chosen metric
(not necessarily related to the connection) is helpful just as in
the theory of affine hypersurfaces or in the geometry of statistical
structures.

\bigskip

\section{Preliminaries}
In this section we  introduce  notions appearing in this paper and
provide basic information on them. All details  can be found in
\cite{NS}, \cite{LSZ} and \cite{BW4}.

 By a statistical structure on
a manifold $M$ we mean a pair $(g, \nabla)$, where $g$ is a positive
definite metric tensor field and $\nabla$ is a torsion-free
connection (called a statistical connection) on $M$ such that the
cubic form $\nabla g$ is symmetric for all arguments. Such a pair
used to be also called a Codazzi pair. Denote by $\hat\nabla$ the
Levi-Civita connection for $g$ and by $K$ the difference tensor
between $\nabla$ and $\hat\nabla$, that is,
\begin{equation}
K_XY=\nabla_XY-\hat\nabla_XY,
\end{equation}
where $X,Y$ are vector fields on $M$. $K(X,Y)$ will stand for
$K_XY$. The $(1,2)$-tensor field $K$ is symmetric and symmetric
relative to $g$. Alternatively, we can use the symmetric cubic form
$A$ defined as  $A(X,Y,Z)=g(K(X,Y),Z))$. One easily sees that
$\nabla g=-2A$. It is  clear that a statistical structure can be
defined as a pair $(g,A)$, where $g$ is a Riemannian metric and $A$
is a symmetric cubic form. A statistical structure is called trivial
if $\nabla=\hnabla$, (equiv. $K\equiv 0$ or $A\equiv 0$).

Of course, one can consider the difference tensor $K=\nabla-\hat
\nabla$ and the cubic forms $\nabla g$ and $A(X,Y,Z)=g(K(X,Y),Z)$
for any connection $\nabla$. If the connection $\nabla$ is
torsion-free then $K$ is symmetric and not necessarily symmetric
relative to $g$. We always have $\nabla
g(X,Y,Z)=-A(X,Y,Z)-A(X,Z,Y)$. In particular, $\nabla
g(X,X,X)=-2A(X,X,X)$ for any $X$.

Let $(g,\nabla)$ be a statistical structure. The dual connection
$\onabla$ for the statistical connection $\nabla$ is defined by the
formula
\begin{equation}
Xg(Y,Z)=g(\nabla_XY,Z)+g(X,\onabla_XZ).
\end{equation}
The dual connection $\onabla$ is torsion-free and  the pair
$(g,\onabla)$ is a statistical connection. The difference tensor for
the dual connection $\onabla$  is equal to $-K$, that is,
$\onabla=\hnabla -K$.

Let $R$, $\oR$, $\hR$ stand for the curvature tensor for $\nabla$,
$\onabla$ and $\hnabla$ respectively. The corresponding Ricci
tensors for $\nabla$ and $\onabla$ will be denoted  by $\Ric$ and
$\overline{\Ric}$. These Ricci tensors do not have to be symmetric.
But they are simultaneously symmetric. Note that $\Ric$ is symmetric
if and only if $d\tau=0$, where $\tau (X)=\tr K_X$.

 The $(0,4)$-tensor field $g(R(X,Y)Z, W)$ is not, in general,
 skew-symmetric for $Z,W$. The class of statistical structures whose
 curvature tensor  fulfills this symmetry condition is very rich, although the condition is very strong.
  Such structures appear in many situations. In particular,
 they exist on affine spheres, which  constitute the most
 distinguished  class of affine hypersurfaces.
 We have
\begin{lemma}\label{przeniesiony_lemat}
Let $(g,\nabla) $ be a statistical structure. The following
conditions are equivalent:
\newline
{\rm 1)} $R=\overline R$,
\newline
{\rm 2)} $\hnabla K$ is symmetric (equiv. $\hat\nabla A$ is
symmetric),
\newline
{\rm 3)} $g(R(X,Y)Z,W)$ is skew-symmetric relative to $Z,W$.
\end{lemma}
  A statistical structure satisfying one of the above conditions
is called conjugate symmetric. Note that the conjugate symmetry
implies the symmetry of $\Ric$. We shall need

\begin{proposition}\label{projectively_flat}
Let $(g,\nabla)$ be a conjugate symmetric statistical structure. The
connection $\nabla$ is projectively flat if and only if $\onabla$ is
projectively flat.
\end{proposition}
\proof Since the statistical structure is conjugate symmetric, the
Ricci tensors for $\nabla$ and $\onabla$ are identical and
symmetric. Consider first the case where $\dim M>2$. A
Ricci-symmetric torsion-free connection is projectively flat if and
only if its curvature tensor is of the form
$$ R(X,Y)Z=\Phi(Y,Z)X-\Phi(X,Z)Y$$
for some symmetric $(0,2)$-tensor field $\Phi$.  Since $R=\overline
R$, the connections $\nabla$ and $\onabla$ are simultaneously
projectively flat. Assume now that $\dim M=2$. In this case the
projective flatness of a connection $\nabla$  is equivalent to the
symmetry of the cubic form $\nabla \Ric$. Assume $\nabla\Ric$ is
symmetric. Since $R$ satisfies {\rm 3)} from Lemma
\ref{przeniesiony_lemat}, on a $2$-dimensional manifold we have
$$R(X,Y)Z=k\{g(Y,Z)X-g(X,Z)Y\},$$
where $k$ is the sectional $\nabla$-curvature of the statistical
structure, see \cite{BW4}. On the other hand, we always have
$$R(X,Y)Z=\Ric (Y,Z)X-\Ric (X,Z)Y.$$
It follows that $\Ric =kg$.  We now have $\onabla
\overline{\Ric}=\onabla \Ric=\nabla\Ric-2K\Ric=\nabla\Ric-2k(K g)$.
Since $K g$ is a symmetric cubic form, we get the desired
condition.\koniec
\bigskip

\medskip
We shall now briefly recall the way in which statistical structures
appear on locally strongly convex affine hypersurfaces.

Let $f :M\to \R ^{n+1}$ be an immersed  hypersurface. For simplicity
assume that $M$ is connected and orientable. Let $\xi$ be a
transversal vector field for the immersion $f$. We define the
induced volume form $\nu _\xi$ on $M$  as follows
$$\nu_\xi (X_1,...,X_n)=\det (f_*X_1,...,f_*X_n,\xi).$$
We  also have the induced connection $\nabla$ (torsion-free) and the
second fundamental form $g$ (symmetric) defined by the Gauss formula
\begin{equation}\label{Gauss_formula}
D_Xf_*Y=f_*\nabla _XY +g(X,Y)\xi,
\end{equation}
 where $D$ is the standard flat
connection on $\R ^{n+1}$. A hypersurface is called non-degenerate
if $g$ is non-degenerate. If the hypersurface is locally strongly
convex, $g$ is definite. By multiplying $\xi$ by $-1$, if necessary,
we can assume that $g$ is positive definite. A transversal vector
field is called equiaffine if $\nabla \nu_\xi=0$. This condition is
equivalent to the symmetry of the cubic form $\nabla g$. In
particular, if the hypersurface is locally strongly convex,
$(g,\nabla)$ is a statistical structure. It follows that for a
statistical structure obtained on a hypersurface by a choice of a
transversal vector field, the Ricci tensor of $\nabla$ is
automatically symmetric. A hypersurface equipped with an equiaffine
transversal vector field is called an equiaffine hypersurface. If
$\xi$ is an equiaffine transversal vector field, the Weingarten
formula looks as follows
\begin{equation}\label{Weingarten-formula}
D_X\xi= -f_*\mathcal SX.
\end{equation}
 The $(1,1)$-tensor field $\mathcal S$ is
called the shape operator for $\xi$.  It is symmetric relative to
$g$. We have the Gauss equation
\begin{equation}\label{Gauss_equation_for R}
R(X,Y)Z=g(Y,Z)\mathcal SX-g(X,Z)\mathcal SY.
\end{equation}
For an equiaffine hypersurface $f$ we define  the conormal map
 $$\overline f: M\to
(\R^{n+1})^*\setminus \{0\}$$ as follows
$$\overline {f}(x)(\xi_x)=1, \ \ \ \ (\overline f(x))_{|f_*(T_xM)}\equiv 0.$$
If $f$ is non-degenerate then $\overline f$ is an immersion. In this
case, for each $x\in M$ the conormal vector $\overline
f_x=\overrightarrow{0\overline{f_x}}$ is transversal to $\overline
f$. We equip the immersion $\overline f$ with the equiaffine
transversal vector field $-\overline f$. Again, we receive the
induced objects on $M$. In particular, the induced connection turns
out to be the dual connection for $\nabla$ relative to $g$. The
second fundamental form for the conormal immersion is equal to
$g(S\cdot,\cdot)$. The Gauss equation for the conormal immersion is
the following
\begin{equation}\label{Gauss_equation_for_oR}
\overline R(X,Y)Z=g(Y,\mathcal SZ)X-g(X,\mathcal SZ)Y.
\end{equation}
It means that the dual connection $\onabla$ is projectively flat if
$n>2$. The dual connection is also projectively flat for $n=2$.

 We
shall use the following version of the fundamental theorem in affine
differential geometry, see e.g. \cite{DNV}.
\begin{thm}\label{fundamental}
Let $(g,\nabla)$ be a statistical structure  on a simply connected
manifold $M$. If $\nabla$ is Ricci-symmetric and the dual connection
$\onabla$ is projectively flat then there is an equiaffine immersion
$f:M\to \R ^{n+1}$ such that $g$ is the second fundamental form and
$\nabla$ is the induced connection for $f$.
\end{thm}

Among equiaffine  hypersurfaces we distinguish equiaffine spheres. A
proper equiaffine sphere is an equiaffine hypersurface for which the
affine lines detrmined by the chosen equiaffine transversal vector
field intersect at one point, called the center of the sphere. The
center does not belong to the hypersurface and the equiaffine
transversal vector field is a non-zero constant multiple of the
position vector field relative to the center.  The shape operator
for a proper equiaffine sphere is a non-zero constant multiple of
the identity.
 A proper equiaffine
sphere is also said to be a centroaffine hypersurface or a
hypersurface  with a centroaffine normalization. As it was noticed
above, if a hypersurface is locally strongly convex, the conormal
map is naturally equipped with a centroaffine normalization. By an
improper equiaffine sphere we mean a hypersurface equipped with a
constant transversal vector field. In this case the shape operator
vanishes. An equiaffine locally strongly convex hypersurface is an
equiaffine sphere (proper or improper) if an only if the induced
statistical structure is conjugate symmetric.

On a locally strongly convex  hypersurface we
 also have the volume form $\nu_g$ determined by $g$. In general,
this volume form is not covariant constant relative to $\nabla$. A
basic theorem of the classical affine differential geometry says
that there is a unique equiaffine transversal vector field $\xi$
such that $\nu_\xi =\nu_g$. This unique transversal vector field is
called the affine normal vector field or the Blaschke affine normal.
The second fundamental form for the affine normal is called the
Blaschke metric or the affine metric. A hypersurface endowed with
the affine Blaschke normal is called a Blaschke hypersurface. An
affine Blaschke sphere (or just an affine sphere) is a Blaschke
hypersurface which is an equiaffine sphere. The easiest examples of
Blaschke affine spheres are ellipsoids, hyperboloids and
paraboloids. More precisely, when we endow an ellipsoid with the
centroaffine normalization with the center at the center of the
ellipsoid and a transversal vector field being the opposite to the
position vector relative to the center, we get a proper Blaschke
affine sphere whose induced statistical structure is trivial.
Hyperboloids can be treated in a similar way. For an elliptic
paraboloid the Blaschke affine normal is parallel to the axis of the
paraboloid (so it is an improper affine sphere) and the induced
statistical structure is also trivial.
 A classical theorem says that the
only compact Blaschke affine spheres are ellipsoids with its trivial
structure.

\bigskip

\section{Geodesics of affine connections}

Before we study geodesics of  statistical connections we shall
collect some facts about geodesics of affine connections in general.
Some of the facts are likely known, but we prove all of them.
Riemannian metrics which we consider in this section are not related
to the connections in any way.

By a $\nabla$-geodesic we mean a geodesic relative to  a connection
$\nabla$ parametrized by its affine parameter, if not otherwise
stated. A geodesic is called complete if its affine parameter runs
from $-\infty$ to $\infty$. A pre-geodesic is a curve, parametrized
or unparametrized, which can be parametrized or reparametrized in
such a way that it becomes a geodesic. Since geodesics are regular
curves, we will consider only regular parametrizations of
pre-geodesics.

 Recall first the following
 theorem proved, for instance, in \cite{KN} in Chapter III, see
Theorem 8.7. For a given connection on $M$  let $(x^1,...,x^n)$ be a
normal coordinate system with origin at $x_0\in M$. Let $\mathcal
U(x_0;\rho)$ denote the geodesic ball of radius $\rho$, that is,
$\mathcal U(x_0;\rho)= \{ (x^1,...,x^n)\in M:
(x^1)^2+...+(x^n)^2<\rho ^2\}$.

\begin{thm}\label{geodesicsKN} Let  $M$ be equipped with a linear connection. For each
point $x_0$ of $M$ and its normal coordinate system $(x^1,..., x^n)$
with origin  at $x_0$ there is a positive number $a$  such that if
$0<\rho<a$ then $\mathcal U(x_0;\rho)$ is geodesically convex and
each point from $\mathcal U(x_0;\rho)$ has a normal coordinate
neighborhood containing $\mathcal U(x_0;\rho)$.
\end{thm}

Using this theorem one can prove

\begin{proposition}\label{cylinder}
Let $M$ be equipped with a connection. If $\gamma:(a,b)\to M$ is a
geodesic and there is the limit of $\gamma(t)$ in $M$ for $t\to b$
then $\gamma$ can be extended as a geodesic beyond
$b$.\end{proposition}

\proof Let $x_0=\lim_{t\to b} \gamma(t)$ and $\mathcal U$ be a
normal neighborhood of $x_0$ as in Theorem \ref{geodesicsKN}. In
this  neighborhood there is  $\gamma (t_0)$ for some $t_0\in (a,b)$.
We can assume that $t_0=0$. There is a normal neighborhood $\mathcal
U '$ of $p:=\gamma(0)$ containing $\mathcal U$. Take the
neighborhood $exp_p^{-1}(\mathcal U')$ of $0$ in $T_pM$. We have
$\beta(t):=exp_p^{-1}(\gamma(t))=\dot{\gamma}(0)t$ for $t\in [0,b)$.
Since $x_0\in \mathcal U\subset\mathcal U'$, the segment $\beta(t)$
of a straight line can be extended beyond $b$. Hence the geodesic
$\gamma$ can be extended beyond $b$.\koniec

By the same arguments as in the last proof we get

\begin{proposition}\label{geodesic_infinity}
Let $M$ be equipped with a connection. If $\gamma:(a,\infty)\to M$
is a geodesic then  there is no limit of $\gamma(t)$ in $M$ for
$t\to\infty$.
\end{proposition}

We  can now reformulate Proposition \ref{cylinder} as follows

\begin{proposition}\label{pregeodesic}
Let $M$ be equipped with a connection. If $\gamma (t)$ for $t\in (a,
b)$, where $b\in \R$ or $b=\infty$ is a pre-geodesic and there is
$x_0=\lim_{t\to b}\gamma(t)$ in $M$ then the pre-geodesic can be
extended as a pre-geodesic beyond $x_0$.
\end{proposition}
\proof Since we can reparametrize the curve to an affine
parameterization, we can assume that $\gamma(t)$ is a geodesic. By
Proposition \ref{geodesic_infinity} we  have that $b\in\R$. It is
now sufficient to Proposition \ref{cylinder}.\koniec
\bigskip

\medskip

Let now $(M,g)$ be a Riemannian manifold. The distance between
points $x,y\in M$ will be denoted by $d(x,y)$ and the length of a
curve $\gamma$ by $d\gamma$.

\begin{proposition}\label{nabla-geodesic}
Let $\gamma:[a,b)\to M$, where $b\in \R$,  be a $\mathcal C^1$-
curve. If the curve is extendable (as a $\mathcal C^1$-curve) to
$[a,b]$ then $\Vert \dot{\gamma} \Vert$ is bounded in $[a,b)$.
 If $(M,g)$ is complete and
 $\Vert
 \dot{\gamma}\Vert$ is bounded on $[a,b)$, then there is the limit $\lim _{t\to
b}\gamma (t)\in M$.
\end{proposition}
\proof The first assertion is trivial. Assume now that $(M,g)$ is
complete and $\Vert \dot{\gamma} \Vert\le N$ in $[a,b)$. Let $t_n\in
[a,b)$ be a sequence converging to $b$. We have

\begin{eqnarray*}
d(\gamma (t_n),\gamma (t_m))&&\le \ the \ length\ of \ \gamma\
between
\ \gamma( t_n) \ and \  \gamma( t_m)\\
&&= \int_{t_n}^{t_m} \Vert \dot{\gamma}(t) \Vert dt\le N|t_m-t_n|.
\end{eqnarray*}
Hence the sequence $\gamma (t_n)$ is a Cauchy one.\koniec

\begin{proposition}\label{well-known} Let $\gamma(t)$ for $t\in [a,b)$, where $b\in \R$ or
$b=\infty$,
 be a curve in a
complete Riemannian manifold $(M,g)$. If $d\gamma<\infty$ then $\lim
\gamma(t)$ for $t\to b$  exists in $M$.
\end{proposition}

 \proof  Let
$d\gamma<N<\infty$. The closed ball $\mathcal B=\{x\in M:\ \ d(x,
\gamma (a))\le N\}$ is compact  and $\im \gamma\subset \mathcal B$.
Take any sequence $t_n\to b$. There is a subsequence $t_{n_m}$ such
that the sequence of points $\gamma(t_{n_m})$ has a limit in
$\mathcal B$. We now claim that $\gamma(t)$ has a limit for $t\to
b$. Suppose that there are two sequences $t_n\to b$ and $s_n\to b$
such that $\gamma(t_n)\to p$ and $\gamma (s_n)\to q$ where $p\ne q$.
We can assume that $t_n$ and $s_n$ are increasing and
$$t_1<s_1<t_2<s_2<t_3... .$$
 Take two balls $\mathcal B_1$, $\mathcal B_2$ with centers at
 $p$ and $q$ respectively, such that $d(\mathcal B_1,\mathcal
 B_2)>\delta>0$. We can assume that $\gamma(t_n)\in \mathcal B_1$
 and $\gamma(s_n)\in \mathcal B_2$ for all $n$. Then
 $d\gamma|^{s_1}_{t_1}>\delta, \ \ d\gamma|^{s_2}_{t_2}>\delta,\ \ ...$. It follows that
 $d\gamma=\infty$, which gives a contradiction.
 \koniec

\begin{proposition}\label{s_to_infty} Let $(M,g)$  be a complete Riemannian manifold and
$\nabla$  a connection on $M$. If $r(s)$ is an arc-length
parametrization of a maximal $\nabla$-geodesic then the parameter
$s$ runs from $-\infty$ to $\infty$. In particular, every maximal
$\nabla$-geodesic has infinite length.
\end{proposition}
\proof Assume that the whole domain of $r$ is $(a,b)$, where $b\in
\R$. By Proposition \ref{nabla-geodesic} we know that there exists
the limit of $r(s)$ for $s\to b$. By Proposition \ref{pregeodesic}
the pre-geodesic can be extended, which is a contradiction with the
maximality. Hence $b=\infty$. For the same reasons $a=-\infty.
$\koniec

\begin{proposition}\label{arc_length}
Let $\nabla$ be a connection on a Riemannian manifold $(M,g)$. If
$r(s)$ is an arc-length parametrization of a  $\nabla$-geodesic
$\gamma(t)$, then
\begin{equation}
\nabla_{\dot r}{\dot r}=A(\dot r,\dot r,\dot r)\dot r.
\end{equation}
Moreover, if $t=\varphi (s)$ is the change of the parameters, where
$\varphi'>0$, then
\begin{equation}
A(\dot r,\dot r,\dot r)=-\frac{d}{ds}\ln \Vert\dot\gamma\circ
\varphi\Vert.
\end{equation}
\end{proposition}
\proof We have
$$ 0=\frac{d}{ds}g(\dot r,\dot r) =\nabla g(\dot r,\dot r, \dot
r)+2g(\nabla_{\dot r}\dot r, \dot r).$$ Since $\nabla _{\dot r} \dot
r$ is parallel to $\dot r$, we now have $\nabla_{\dot r}\dot
r=-\frac{1}{2}\nabla g(\dot r,\dot r,\dot r)\dot r=A(\dot r,\dot
r,\dot r)\dot r$.

For proving the second assertion observe that since
$r(s)=\gamma(\varphi(s))$ and $\dot
r(s)=\varphi'(s)(\dot\gamma\circ\varphi)(s)$, we have
\begin{equation}\label{varphi'}
\Vert(\dot\gamma\circ\varphi)(s)\Vert=\frac{1}{\varphi'(s)}=g((\dot\gamma\circ\varphi)(s),\dot
r(s)).
\end{equation}
Consequently
\begin{eqnarray*}
&&\Vert\dot\gamma\circ\varphi\Vert'=\frac{d}{ds}g(\dot\gamma\circ\varphi, \dot r)\\
 &&\ \ \ \ \ \ \ \ \ \  = \nabla g(\dot r,\dot\gamma\circ \varphi ,\dot r)) +
 g(\nabla_{\dot r}(\dot\gamma\circ\varphi),\dot
 r)+g(\dot\gamma\circ\varphi,\nabla_{\dot r}{\dot r})\\
 &&\ \ \ \ \ \ \ \ \ \ =\frac{1}{\varphi '}\nabla g(\dot r,\dot r, \dot r)+g(
 \varphi '\nabla_{\dot\gamma\circ\varphi}(\dot\gamma\circ\varphi),\dot
 r)+\frac{1}{\varphi '}A(\dot r,\dot r,\dot r)\\
 &&\ \ \ \ \ \ \ \ \ \  =-2\frac{1}{\varphi '}A(\dot r,\dot r,\dot r)  +\frac{1}{\varphi '}A(\dot r,\dot r,\dot
 r)=-\frac{1}{\varphi '}A(\dot r,\dot r,\dot r).
\end{eqnarray*}
Take $\Lambda=-\ln \Vert\dot\gamma\circ \varphi\Vert$. Then
$$\frac{d}{ds}\Lambda =-\frac{\Vert\dot\gamma\circ \varphi\Vert'}{\Vert\dot\gamma\circ
\varphi\Vert}=\frac{\frac{1}{\varphi '}A(\dot r,\dot r,\dot
r)}{\frac{1}{\varphi'}}.$$

\koniec

\begin{proposition}
Let  $r(s)$, $s\in\R$, be any regular parametrization of a maximal
geodesic of some affine connection $\nabla$ on a manifold $M$. Let
$\nabla_{\dot r}{\dot r}=\Lambda' \dot r$, where $\Lambda $ is a
function  bounded from below or above on the whole of $\R$. Then the
geodesic is complete.
\end{proposition}
\proof Let $t=\varphi(s)$ be an affine parameter of our geodesic,
say $\gamma (t)$. We can assume that $\varphi'>0$ on $\R$.
Then $\gamma(\varphi(s))=r(s)$, $\varphi '(\dot\gamma
\circ\varphi)=\dot r$ and consequently
\begin{equation}
\nabla_{\dot r}\dot r=\nabla_{\dot
r}(\varphi'(\dot\gamma\circ\varphi))=\varphi''(\dot\gamma\circ
\varphi).
\end{equation}
Hence $\Lambda'=(\ln \varphi')' $ and
 we can assume that
$ \varphi'=e^\Lambda . $ If $\Lambda$ is bounded from below then
$\varphi'$ is greater than some positive number. Therefore
$\varphi(s)_{s\to\infty}\to\infty$ and
$\varphi(s)_{s\to-\infty}\to-\infty$. If $\Lambda $ is bounded from
above then we can change the orientation of the parametrization
$r(s)$ and in the equation $\nabla_{\dot r}{\dot r}=\Lambda' \dot r$
we replace $\Lambda$ by $-\Lambda$, which is now bounded from below.
\koniec
\medskip

\begin{proposition} Let $\nabla$ be a connection on a complete Riemannian
manifold. If a maximal $\nabla$-geodesic has scalar speed bounded
(from above), then the geodesic is complete.
\end{proposition}

\proof By Proposition \ref{s_to_infty} we know that the arc-length
parameter of the geodesic runs from $-\infty$ to $+\infty$. Use the
same notations and  agreements as in Proposition \ref{arc_length}.
By (\ref{varphi'}) $\varphi'$ is greater than some positive number
for every $s$. It means that $\varphi(s)_{s\to\infty}\to\infty$ and
$\varphi(s)_{s\to-\infty}\to-\infty$. \koniec .

For a $\nabla$-geodesic $\gamma(t)$ on a Riemannian manifold we
introduce two functions $l(t)=\Vert\dot\gamma(t)\Vert$ and
$u(t)=\frac{\dot\gamma(t)}{\Vert\dot\gamma(t)\Vert}$.
 We have
\begin{equation}\label{l^2'}
\begin{array}{rcl}&&\frac{d}{dt}g(\dot\gamma (t),\dot\gamma(t)) =
2g(\hat\nabla_{\dot\gamma(t)}
\dot\gamma,\dot\gamma(t))=-2g(K_{\dot\gamma
(t)}\dot\gamma(t),\dot\gamma(t))\\
&&=-2A(\dot\gamma (t), \dot\gamma (t),\dot\gamma
(t))=-2A\left(\frac{\dot\gamma(t)}{\Vert\dot\gamma(t)\Vert},
\frac{\dot\gamma(t)}{\Vert\dot\gamma(t)\Vert},
\frac{\dot\gamma(t)}{\Vert\dot\gamma(t)\Vert}\right)\Vert\dot\gamma(t)\Vert
^3.
\end{array}
\end{equation}
Hence
\begin{equation}\frac{d}{dt}(l^2(t))=-2A(u(t),u(t),u(t))l^3(t)\end{equation}
 and therefore
\begin{equation}\label{first_derivative1/l}
\left(\frac{1}{l}\right)'=A(u,u,u).
\end{equation}

Assume that the geodesic $\gamma$ is defined on the interval
$[a,b)$. We   can assume that  $a=0$ and $l(0)=\Vert\dot{\gamma
}(0)\Vert =1$. Assume that
\begin{equation}\label{-N_N}
-N\le A(U,U,U)\le N\end{equation}
 for some
non-negative number $N$ and all $U\in \mathcal UM$, where $\mathcal
UM$ is the unit sphere bundle over $M$. It happens, for instance, if
$M$ is compact. We have

\begin{equation}
\frac{1}{l(t)}=1+A(u(\theta),u(\theta), u(\theta))t
\end{equation}
for every  $t\in (0,b)$ and some $\theta\in (0, t)$.
 Since
$A(u(\theta),u(\theta), u(\theta)) \ge -N$, we get
\begin{equation}
l(t)\le \frac{1}{1-N t},
\end{equation}
if $t< \frac{1}{N}$. The function $[0, b)\ni t\to \frac{1}{1-N t}$
is positive valued and increasing. Using now Proposition
\ref{nabla-geodesic} we get

\begin{proposition}
Let $\nabla$ be a connection on a Riemannian manifold $M$ whose
metric  is complete and {\rm (\ref{-N_N})} holds. Let
$\gamma:[0,b)\to M$ be a $\nabla$-geodesic such that $\Vert
\dot{\gamma}(0)\Vert=1$. Then the geodesic can be extended at least
up to any  parameter $c<\frac{1}{N}$.
\end{proposition}
\proof Indeed, if $c<\frac{1}{N}$, then the function $l$ is bounded
on the interval $[0, c)$. \koniec

In contrast with Proposition \ref{s_to_infty} in the following
proposition the Riemannian manifold does not have to be complete.
\begin{proposition}
Let $\nabla$ be a connection  on a Riemannian manifold and
{\rm(\ref{-N_N})} holds. Every complete $\nabla$-geodesic has
infinite length.
\end{proposition}
\proof By (\ref{first_derivative1/l}) we have
$\left(\frac{1}{l}\right)'\le N$. By integrating this inequality on
an interval $[0, t]$ we get $\frac{1}{l(t)}\le Nt +1$. Consequently,
for every $t\in R^+$ we have $ l(t)\ge \frac{1}{Nt+1}.$ Integrating
this inequality on an interval $[0,t]$ we obtain
\begin{equation}
d\gamma_{[0,t]}\ge \frac{\ln(Nt+1)}{N}.
\end{equation}

\koniec

\medskip

It is easy to find  an incomplete metric and its statistical
connection whose complete geodesic has finite length.

\begin{example}{\rm Let $M$ be the open strip $ \R\times(-L,L) \subset \R^2$
endowed with the standard metric $g$, where
$L=\frac{\sqrt{\pi}}{2}$. Let $G(t)$ be the antiderivative of the
Gaussian function $e^{-t^2}$ with $G(0)=0$. Denote by $(U,V)$ the
canonical frame on $M$. Each point of $M$ has the coordinates
$(x,G(t))$ for some $t\in \R$. Take any symmetric and symmetric
relative to $g$ $(1,1)$-tensor field   $K$  such that
\begin{equation}
K_{(x,G(t)) }(V,V)=2te^{t^2}V.
\end{equation}
Let $\nabla=\hat\nabla +K$. The piece of a straight line
$\gamma(t)=(0, G(t))$ for $t\in (-\infty,\infty)$ is a complete
$\nabla$-geodesic. Of course, it has finite length $\sqrt{\pi}$.
Note that {\rm(\ref{-N_N})} is not satisfied.}
\end{example}

 We shall now give an example of a statistical structure on a compact
manifold whose statistical connection is not complete. The example
is important for the last  section of this paper.
\begin{example}\label{example_compact_non-complete}
{\rm Take  $\R^2$ with its standard flat Riemannian structure. Let
$U,V$ be the canonical frame field on $\R^2$. Define the statistical
connection $\nabla$ as follows

\begin{equation}
\nabla_UU=U, \ \ \nabla_UV=-V , \ \  \nabla _VV=-U.
\end{equation}
 The statistical structure can be projected on the standard torus
$T^2=\R^2/\Z ^2$.  Observe that $\hat\nabla K=0$ (equivalently
$\hat\nabla A=0$), hence the structure is conjugate symmetric.  For
later use, observe also that $\nabla$ is projectively flat and its
Ricci tensor is non-degenerate negative definite.

A curve $\gamma(t)=(x(t), y(t))$ is a $\nabla$-geodesic if and only
if

\begin{equation}
x''+(x')^2-(y')^2=0,\ \ \ y''-2x'y'=0.
\end{equation}
Let  $y_0$ be a fixed real number. Consider the curve
\begin{equation}
\gamma(t)=(\ln(1-t), y_0)
\end{equation}
for $t\in [0, 1)$. It is a $\nabla$-geodesic. We have $\gamma(t)\to
(-\infty, y_0)$, $\Vert \dot{\gamma }(t)\Vert=\frac{1}{1-t}\to
+\infty$ for $t\to 1$. This geodesic cannot be extended beyond $1$
and its image   is closed in $T^2$. It means, in particular, that
closed geodesics  of affine connections do not have to be complete.
}\end{example}

\bigskip

\section{Completeness of statistical connections}

 The following theorem generalizes Example
\ref{example_compact_non-complete}.
\begin{thm}\label{non-complete}
Let $(g,\nabla)$ be a statistical structure for which $\nabla$ is
complete and $(\hat\nabla A)(U,U,U,U)\le 0$ for each $U\in\mathcal
UM$. Then the statistical structure must be trivial.\end{thm} \proof
We shall first prove
\begin{lemma}\label{second_derivative_1/l}
Let $\gamma(t)$  be a $\nabla$-geodesic. Then

\begin{equation}\label{derivative1/l}
\begin{array}{rcl}
&& (A(u,u,u))'(t)\\
&&\ \ \ \ \ \ \ \ =l(t)[(\hat\nabla
A)(u(t),u(t),u(t),u(t))-3\sum_{i= 2}^n(A(u(t),u(t),e_i(t)))^2],
\end{array}
\end{equation}
where $e_1(t)=u(t), e_2(t),..., e_n(t)$ is an orthonormal basis of
$T_{\gamma(t)}M$. In particular,

\begin{equation}
(A(u,u,u))'(t)\le l(t)(\hat\nabla A)(u(t),u(t),u(t),u(t)).
\end{equation}
\end{lemma}
\proof Using (\ref{first_derivative1/l}) we get

\begin{equation}
\begin{array}{rcl}
&&\hat\nabla_{\dot\gamma(t)}u=\hat\nabla_{\dot\gamma(t)}\frac{\dot\gamma}{l}\\
&&\ \ \ \  =\left(\frac{1}{l}\right)'(t)\dot\gamma
(t)-\frac{1}{l(t)}K(\dot\gamma(t),\dot\gamma(t))\\
&&\ \ \ \  =l(t)[A(u(t),u(t),u(t))u(t)-K(u(t),u(t))].
\end{array}
\end{equation}
Since
$$K(u(t),u(t))=\sum_{i=1}^ng(K(u(t),u(t)),e_i(t))e_i(t)=\sum_{i=1}^nA(u(t),u(t),e_i(t))e_i(t)$$
for any orthonormal frame  $e_1(t), ..., e_n(t)$,  we now obtain
(using also the symmetry of $A$)
\begin{eqnarray*}
&&(A(u,u,u))'(t)=\dot\gamma (t)(A(u,u,u))\\
&&\ \ \ =(\hat \nabla
A)(\dot\gamma(t),u(t),u(t),u(t))+3A(\hat\nabla_{\dot\gamma(t)}u,u(t),u(t))\\
&&\ \ \ \ =l(t)[(\hat\nabla A)(u(t),u(t),u(t),u(t))
+3(A(u(t),u(t),u(t)))^2\\&&\ \ \ \ \ \ \ \ \ \ \ \ \ \ \ \ \ \
-3(A(u(t),u(t),e_1(t)))^2-3\sum_{i=2}^n(A(u(t),u(t),e_i(t)))^2].
\end{eqnarray*}
By setting $e_1(t)=u(t)$ we complete the proof of the lemma.\koniec

\medskip

 Suppose now  that the structure $(g,\nabla)$ is not trivial. Let $\gamma
:[0,\infty)$ be a $\nabla$-geodesic such that $A(u(0),
u(0),u(0))=-N$, where $N$ is a positive number. The function
$\frac{1}{l}$ is smooth and positive on $[0,\infty)$. By Taylor's
formula and (\ref{first_derivative1/l}) we have
\begin{equation}\label{Taylor_second}
\frac{1}{l(t)}=\frac{1}{l(0)} -Nt+\alpha(\theta)t^2
\end{equation}
for any $t\in[0,\infty)$, where $\theta\in (0,t)$  and
$2\alpha=(A(u,u,u))'$. Take $t_0>\frac{1}{l(0)N}$. For $t=t_0$ the
expression on the right hand side of (\ref{Taylor_second}) is
negative by Lemma \ref{second_derivative_1/l}. This gives a
contradiction. \koniec

 Observe that on a compact $M$ if $(\hat\nabla A)(U,U,U,U)\le 0$
for every $U\in\mathcal UM$, then $\hat\nabla A= 0$ on $M$. It
follows from Ros' integral formula. Namely,   the formula says that
for any covariant tensor field $s$ on a compact Riemannian manifold
we have $\int_{\mathcal{U}M} (\hat\nabla s)(U,...,U)dU=0.$ However,
the fact that $\hat\nabla A=0$ does not trivialize the situation.
First of all, even in the theory of affine hypersurfaces one knows
examples of non-trivial statistical structures, for which
$\hat\nabla A=0$. The most famous is that on the hypersurface of
$\R^{n+1}$ given by the equation $x_1\cdot ...\cdot x_{n+1}=1$ for
$x_1>0,...,x_{n+1}>0.$ This hypersurface is locally strongly convex
affine sphere with $\hat\nabla A=0$, see \cite{LSZ}.
 By Theorem \ref{non-complete}
neither the induced connection $\nabla$ nor its dual are complete.
On the other hand, the Blaschke metric on this hypersurface is
complete, see \cite{LSZ}. The class of Blaschke  hypersurfaces
satisfying the condition $\hat\nabla A=0$ turns out to be especially
interesting. In the paper \cite{HLV} all such hypersurfaces were
classified. By Theorem \ref{non-complete} we know that  in all
non-trivial  cases  the induced connection is not complete.

 Note that  the compact case of hypersurfaces satisfying the condition $\hat\nabla A=0$
 is not interesting. Namely, for an ovaloid this
condition implies that the ovaloid must be an ellipsoid with its
trivial affine structure. More precisely, we have
\begin{proposition}
Let $:M\to \R^{n+1}$ be an equiaffine ovaloid  with the induced
statistical structure $(g,\nabla)$. If $\hat\nabla A=0$ then the
structure is trivial and the ovaloid is an
ellipsoid.\end{proposition}

\proof Since $\hat\nabla A$ is symmetric, the structure is conjugate
symmetric, and consequently, the hypersurface is an equiaffine
sphere and the Ricci tensor $\Ric$ is symmetric. It follows that
$d\tau=0$, where $\tau (X)=\tr K_X$. Hence there is a function
$\rho$ on $M$ such that $\tau =d\rho$. There is a point $p\in M$ at
which $\rho$ attains an extremum. Then $\tau_p=0$. Since $\hat\nabla
K=0$, we have  $\hat\nabla \tau=0$ and consequently $\tau=0$ on $M$.
It means that the hypersurface is a compact Blaschke affine sphere.
It is well-known from the classical affine differential geometry
that such a sphere must be an ellipsoid with its trivial
structure.\koniec

\medskip

We shall now prove a positive result, that is, a theorem showing how
to produce complete affine connections.

\begin{thm}\label{Noguchi-general}
 Let $(M,g)$ be a complete Riemannian manifold and $A$  a
cubic form given by
\begin{equation}
A(X,Y,Z)= \alpha _1g(X,Y)d\sigma(Z)+\alpha_2 g(Y,Z)d\sigma(X)
+\alpha_3 g(X,Z)d\sigma(Y)
\end{equation}
 for some
function $\sigma$ on $M$ and real numbers $\alpha_1$, $\alpha _2$,
$\alpha_3$ such that $\alpha=\alpha_1+\alpha _2+\alpha_3\ge 0$.
Assume that the function $\sigma$ is bounded from below on $M$. Let
$K$ be a $(1,1)$-tensor field defined by the formula $g(K(X,Y),
Z)=A(X,Y,Z)$. Then the connection $\nabla =\hat\nabla+K$  is
complete.
\end{thm}
\proof Let $\gamma:[0,b)\to M$ be a $\nabla$-geodesic such that
$l(0)=1$. We have $A(\dot{\gamma},\dot{\gamma},\dot{\gamma})=\alpha
d\sigma(\dot{\gamma})g(\dot{\gamma},\dot{\gamma})
=\alpha{(\sigma\circ \gamma)'}l^2$. Set $L=l^2$. We proved in
(\ref{l^2'}) that $L'= -2A(\dot{\gamma},
\dot{\gamma},\dot{\gamma}).$ Hence $(\ln
L)'=-2\alpha(\sigma\circ\gamma)'$, that is, $
L(t)=e^{-2\alpha\sigma(\gamma(t))+2\alpha\sigma(\gamma(0))}.$ If
$\sigma\ge N$ on $M$ then $L(t)\le e^{2\sigma(\gamma(0))}e^{-2\alpha
N}$. We can now apply Propositions \ref{nabla-geodesic} and
\ref{cylinder}. \koniec

In particular, we have Noguchi's theorem from \cite{No}.

\begin{corollary}\label{Noguchi}
 Let $(M,g)$ be a complete Riemannian manifold and $A$  a
cubic form given by $A=sym (d\sigma\otimes g) $ for some function
$\sigma$ on $M$. Assume that the function $\sigma$ is bounded from
below on $M$. Then the statistical connection of the statistical
structure $(g,A)$ is complete.
\end{corollary}

Note that if the function $\sigma$ is bounded from above then
$-\sigma$ is bounded from below and, consequently, the dual
connection $\onabla$ is complete. In particular, we have

\begin{corollary} Let $(M,g)$ be a compact Riemannian manifold. Each
function $\sigma$ on $M$ gives rise to a statistical structure whose
statistical connection and its dual are complete. The cubic form of
the structure is given by $A=sym(d\sigma\otimes g)$.
\end{corollary}

\medskip

\section{Completeness of the induced connections on equiaffine hypersurfaces}

It is known that pre-geodesics of  induced connections on Blaschke
affine spheres are obtained
 by intersecting  the spheres with affine planes incident with the
 center of the sphere, see \cite{NS}. A more general statement holds.
 Namely, we have

\begin{proposition}\label{pregeodesics} Let $f:M\to\R^{n+1}$ be a
hypersurface equipped with a centroaffine normalization or with a
constant transversal vector field.  The image by $f$ of each
geodesic of the induced connection is contained in a certain affine
plane containing  the center of the centroaffine normalization or,
respectively, is parallel to the constant transversal vector field.
Conversely, if the image by $f$ of a curve in $M$ lies in some
affine plane  containing  the center of the centroaffine
normalization or, respectively, is parallel to the constant
transversal vector field, then the curve is a pre-geodesic of the
induced connection.
\end{proposition}
\proof Denote by $\xi$ the transversal vector field.  Let $\gamma
(t)$ be a geodesic relative to the induced connection $\nabla$. By
the Gauss formula (\ref{Gauss_formula}) we have $ D_{\dot\gamma}
f_*\dot{\gamma}=g(\dot{\gamma}, \dot{\gamma})\xi_\gamma .$ The
Weingarten formula says that $D_{\dot\gamma }\xi=\lambda
f_*\dot\gamma$ for some real number $\lambda$. Denote by $\pi_t$ the
affine plane passing through $f(\gamma(t))$ and whose direction is
spanned by the vectors $f_*\dot{\gamma}(t)$ and $\xi_{\gamma(t)}$.
The plane field $\pi_t$ is parallel relative to the connection $D$
and locally contains $f\circ\gamma$, hence it is constant along
$\gamma$. The entire curve $f\circ\gamma$ lies in the plane.

Conversely, assume that for a curve $\gamma(t)$  its image by $f$
lies in a ceratin affine plane containing the center of the
centroaffine normalization or, respectively, being parallel to the
constant transversal vector field. Of course, $(f\circ \gamma)' (t)$
and $(f\circ\gamma)''(t)$ must lie in the direction of the plane for
every $t$. Such a plane  also contains all transversal vectors
$\xi_{\gamma(t)}$. The direction of the plane is spanned by
$(f\circ\gamma )'(t)$ and $\xi_{\gamma(t)}$ for every $t$.
Therefore, by the Gauss formula, we have $\nabla_{\dot
\gamma}\dot\gamma=\varphi \dot\gamma$ for some function $\varphi$
along $\gamma$. Thus $\gamma(t)$ is a $\nabla$-pre-geodesic.\koniec

Note that a hypersurface from the above theorem does not have to be
non-degenerate. The following lemma will  needed.

\begin{lemma}\label{inR^2}
Let $r(t)$ for $t\in \R$ be  a curve in the vector space  $\R^2$
such that the function $\det (r(t), \dot{r}(t))\ne 0$  is bounded
from $0$ by a non-zero number on the whole of $\R$. There is a
reparametrization of the curve being a diffeomorphism of $\R$ and
 such that after the reparametrization $\ddot{r}$ is parallel to $r$ everywhere on $\R$.
\end{lemma}

\proof  By replacing the parameter $t$ by $-t$, if necessary, we can
assume that $\det (r(t), \dot{r}(t))>0$ and this function is bounded
from $0$ by some positive number for all $t\in\R$. We can write
\begin{equation}
\ddot{r}(t)=\alpha (t) \dot{r}(t)+\beta(t)r(t).
\end{equation}
Set $\rho(t):=\ln \det (r(t),\dot{r}(t))$. One  sees that
\begin{equation}
\rho'=\frac{\det(r,\ddot{r})}{\det(r,\dot{r})}=\alpha.
\end{equation}
We are looking for a reparametrization $s(t)$ such that
$(r(t(s)))''$ is parallel to $r(t(s))$ for every $s$. We have
\begin{equation}
\begin{array}{rcl}
&&(r(t(s)))'=\frac{dt}{ds}\dot{r}(t(s))\\
&&(r(t(s)))''=(\frac{dt}{ds})^2\ddot{r}(t(s))+\frac{d^2t}{ds^2}\dot{r}(t(s))\\
&&\ \ \ \ \ \ \ \ \ =(\frac{dt}{ds})^2
[\alpha(t(s))\dot{r}(t(s))+\beta(t(s))r(t(s))]+\frac{d^2t}{ds^2}\dot{r}(t(s)).
\end{array}
\end{equation}
We want that \begin{equation}\alpha
(t(s))\left(\frac{dt}{ds}\right)^2 +\frac{d^2t}{ds^2}=0.
\end{equation}
The last equality can be written as
\begin{equation}
(\rho(t(s))'\frac{dt}{ds}+\frac{d^2t}{ds^2}=0.
\end{equation}
Denote the function $t(s)$ by $\varphi$ for a moment. Assume that
$\varphi'>0$ everywhere. We have $
(\rho\circ\varphi)'=-\frac{\varphi''}{\varphi'}$ and hence $ \rho
\circ\varphi =-(\ln \varphi')+C,$
 where $C$ is a
constant. Hence $ \varphi'=e^{-\rho\circ \varphi}e^C$ and
consequently
\begin{equation}\label{dspodt}
\frac{ds}{dt}=e^{\rho} e^{-C}=\det (r,\dot{r})e^{-C}. \end{equation}
In order to get a desired reparametrization it is now sufficient to
integrate the function $\det (r,\dot{r})e^{-C}$. The antiderivative
is strictly increasing. Since $\det(r,\dot r)$ is bounded from zero
by some positive number, using (\ref{dspodt}) we get $
\frac{ds}{dt}>N$ for some positive number $N$. It follows that
$s(t)\to\infty$ for $t\to \infty$ and $s(t)\to-\infty$ if $t\to
-\infty$, that is, the reparametrization is a global diffeomorphism
of $\R$.\koniec

 When we have a hypersurface of $\R^{n+1}$ with a
centroaffine normalization then we can assume that the center of the
normalization is $0\in \R^{n+1}$. An affine plane passing through
$0$ is a vector plane and we can assume that it is endowed with some
arbitrary chosen determinant $\det$. We can also assume that the
transversal vector field is equal to the minus position vector field
relative to the center. We have

\begin{thm}\label{procedure}
If for a given hypersurface $f:M\to\R^{n+1}$ endowed with a
centroaffine normalization with  center at $0\in \R^{n+1}$
 each curve $f\circ
\gamma$ lying in a plane $\pi$ containing $0$
 can
be parametrized globally on $\R$ as $\gamma(t)$ in such a way that
the function $\det ((f\circ \gamma)(t), {(f\circ \gamma)'}(t))$ is
bounded from zero by a non-zero number for all $t\in \R$, then the
induced connection $\nabla$ on $M$ is complete.
\end{thm}
\proof  By Proposition \ref{pregeodesics} we know that for every
$\nabla$-pre-geodesic its image by $f$ lies in a certain plane
$\pi\ni 0$. We can parametrize this curve in $\pi$ as in Lemma
\ref{inR^2}. Assume that $f\circ\gamma $ is  such a global
parametrization. We can write $D_{\dot\gamma}f_*\dot\gamma =(f\circ
\gamma)'' =\rho(f\circ \gamma)$ for some function $\rho$ along
$\gamma$. On the other hand
$D_{\dot\gamma}f_*\dot\gamma=f_*\nabla_{\dot\gamma}{\dot\gamma}-g(\dot\gamma,\dot\gamma)(f\circ\gamma)$.
It follows that $\nabla_{\dot\gamma}\dot\gamma=0$.\koniec

We can now prove

\begin{thm}\label{ovaloids}
On a centroaffine ovaloid in $\R^{n+1}$ the induced connection and
its dual  are complete.
\end{thm}
\proof An  ovaloid is a non-degenerate  hypersurface modeled on a
Euclidean sphere. But it is automatically a strongly convex
imbedding of a Euclidean sphere.  Hence we can assume that it is, in
particular, a subset of $\R^{n+1}$. A pre-geodesic of the induced
connection is a curve diffeomorphic to the ordinary circle lying in
some plane $\pi$. Take a periodic regular parametrization $\gamma
(t)$ of a pre-geodesic defined on the whole of $\R$.  Then the
function $\det(\gamma,\dot{\gamma})$ is nowhere zero. Since it is
periodic, it is bounded from $0$ by a non-zero number. We can now
use Theorem \ref{procedure}. Since the affine shape operator for a
centroaffine ovaloid is a non-zero multiple of the identity, the
conormal hypersurface is also locally strongly convex. Hence it is
an ovaloid with the standard centroaffine normalization. The dual
connection is the induced connection for the centroaffine conormal
ovaloid. \koniec

Instead of an ovaloid one can consider a  convex (not necessarily
strongly) imbedding of a Euclidean sphere.  When we have its
centroaffine normalization  then we see, as in the above theorem,
that the induced connection is complete.

When we have an ovaloid in $\R^{n+1}$ equipped with any equiaffine
transversal vector field then the conormal map is an immersion but
not necessarily non-degenerate, that is,  the conormal, in general,
is not an ovaloid in the dual space. But if the shape operator $S$
for the given ovaloid is non-singular everywhere, then the conormal
is a non-degenerate immersion and therefore it is a centroaffine
ovaloid. Hence, by Theorem \ref{ovaloids}, we obtain

\begin{thm}\label{ovaloid2}
For an ovaloid  in $\R^{n+1}$ equipped with any equiaffine
transversal vector field, whose   affine Gauss curvature $\det S$ is
nowhere zero, the dual connection  is complete.
\end{thm}

Lemma \ref{procedure} can be also applied to non-compact manifolds
as well.

\begin{example}{\rm
Consider an elliptic paraboloid in $\R^3$. After choosing a suitable
coordinate system $(x,y,z)$ the paraboloid can be viewed as the
graph of the function $z=x^2+y^2-1$. Consider the centroaffine
normalization with center at $0\in \R^{n+1}$. If we intersect the
paraboloid by a plane passing through the center, we obtain either
an ellipse or a parabola. The ellipses can be treated as in the
proof of Theorem \ref{ovaloids}. Considering  the parabolas, it is
sufficient to look at the parabola $z=x^2-1$ on the plane $y=0$. Its
easiest parametrization is $\gamma(t)=(t, 0,t^2-1)$ for $t\in R$. We
have $\det (\gamma(t),\dot \gamma(t))=t^2+1$. We can now apply
Theorem \ref{procedure}. }\end{example}

Using Theorems \ref{ovaloids} and   \ref{fundamental} one gets

\begin{thm}\label{proj_flat}
Let $(g,\nabla)$ be a statistical structure on a manifold $M$
diffeomorphic to a Euclidean sphere.
\newline
\it 1) If the Ricci tensor of $\nabla$ is symmetric non-degenerate
and $\nabla$ is projectively flat then $\nabla$ is complete on $M$.
\newline
2) If the structure is conjugate symmetric and the  connection
$\nabla$ is projectively flat then $\nabla$ and its dual connection
$\overline\nabla$ are complete on $M$.
\end{thm}
\proof {\it 1)} Since $\Ric$ is symmetric, so is $\overline{\Ric}$.
The dual connection for $\onabla$ is equal to $\nabla$ and this
connection is projectively flat. We can now apply Theorem
\ref{fundamental} for $(g,\onabla)$. We get an ovaloid $f:M\to \R
^{n+1}$ with some equiaffine transversal vector field, the induced
connection $\onabla$ and the second fundamental form $g$. Take the
conormal map $\overline f$ equipped with the standard centroaffine
normalization. Then $\nabla$ is the induced connection for this
centroaffine conormal immersion. Let $h$ be the second fundamental
form for the centroaffine conormal immersion. The shape operator for
the centroaffine conormal immersion is $\id$. The Gauss equation
(for the centroaffine conormal) is the following
\begin{equation}
R(X,Y)Z=h(Y,Z)X-h(X,Z)Y.
\end{equation}
The projective flatness of $\nabla$ implies
\begin{equation}R(X,Y)Z=\frac{1}{n-1}(\Ric(Y,Z)X -\Ric(X,Z)Y)\end{equation}
for $n>2$. For $n=2$ the above equality always holds. Hance
$h=\frac{1}{n-1}\Ric$  is non-degenerate and therefore the
centroaffine conormal immersion is a centroaffine ovaloid with the
induced connection $\nabla$. We can now apply Theorem
\ref{ovaloids}.

{\it 2)} For a conjugate symmetric statistical structures both
statistical connections $\nabla$, $\onabla$ are Ricci-symmetric. By
Proposition \ref{projectively_flat}   both connections  are
projectively flat. By Theorem \ref{fundamental}  the statistical
structure $(g,\nabla)$ can be realized on an ovaloid. Since the
statistical structure is conjugate symmetric, the ovaloid is an
equiaffine sphere. It means that it is a centroaffine ovaloid and,
by Theorem \ref{ovaloids}, $\nabla$ and $\overline\nabla$ are
complete. \koniec

As concerns {\it 1)} in the above theorem note that, by Theorem
\ref{ovaloids}, the connection dual to $\nabla$ relative to
$h=\frac{1}{n-1}\Ric$ is also complete, but it is not the connection
$\overline\nabla$ in general. In the context of Theorem
\ref{proj_flat} let us recall the structure on the torus from
Example \ref{example_compact_non-complete}. That structure  is
conjugate symmetric, the statistical connection is projectively
flat, its Ricci tensor is non-degenerate and the connection is not
complete. The dual connection is not complete either.


\begin{thebibliography}{20}

\bibitem{DNV}  Dillen F., Nomizu K., Vrancken L., \emph{Conjugate connections and
Radon's theorem in affine differential geometry}, Monat. Math., 109
(1990), 221-235


\bibitem{C} Calabi E., \emph{Complete affine hypersurfaces I}, Symposia Math. 10, (1972)
19-38

\bibitem{CY} Cheng S. Y., Yau S.T., \emph{Complete affine hypersurfaces,
Part I, The completeness of affine metrics,} Commun. Pure and Appl.
Math., 39 (1986) 839-866


\bibitem{HLV} Hu Z., Li H., Vrancken L.\emph{ Locally strongly convex affine hypersurfaces with
parallel cubic form}, J. Differential Geometry, 87 (2011) 239-307

\bibitem{KN}  Kobayashi S., Nomizu K., \emph{Foundations of
Differential Geometry} vol. I, John Wiley and Sons, New York, 1963

\bibitem{LSZ} Li A.-M., Simon U., Zhao G.,  \emph{Global Affine Differential Geometry of Hypersurfaces},
W. de Greuter, Berlin-New York, 1993

\bibitem{NS}Nomizu K., Sasaki T., \emph{Affine Differential
Geometry, Geometry of Affine Immersions}, Cambridge University
Press, 1994

\bibitem{No} Noguchi M., \emph{Geometry of statistical manifolds},
Diff. Geom. Appl., 2 (1992), 197-222

\bibitem{BW4} Opozda B., \emph{Bochner's technique for statistical
structures}, Ann. Glob. Anal. Geom., 48 (2015) 357-395

\bibitem{calabi} Opozda B., \emph{Curvature bounded conjugate symmetric
statistical structures with complete
 metric}, Ann. Glob. Anal.
Geom., 55 (2019) 687-702


\end{thebibliography}
\end{document}